\documentclass[12pt,reqno]{amsart}

\usepackage{amssymb}
\usepackage{pdfpages}

\newcommand{\ds}{\displaystyle}     

\DeclareMathAlphabet{\E}{U}{eus}{m}{n}     

\newcommand{\V}{{\mathcal V}}

\newcommand{\PP}{{\mathbb P}}
\newcommand{\C}{{\mathbb C}}
\newcommand{\N}{{\mathbb N}}

\newcommand{\kk}{{\Bbbk}}

\renewcommand{\char}{\mbox{\rm char}}

\newcommand{\Q}{{\mathfrak Q}}

\renewcommand{\sl}{{\mathfrak s}{\mathfrak l}}

\newcommand{\la}{\langle}
\newcommand{\ra}{\rangle}

\newtheorem{thm}{Theorem}[section]

\newtheorem{prop}[thm]{Proposition}

\theoremstyle{definition}
\newtheorem{defn}[thm]{Definition}

\newtheorem{example}[thm]{Example}
\newtheorem{examples}[thm]{Examples}

\newtheorem{ack}{Acknowledgments\!\!\!}  

\newcommand{\dual}{\raisebox{0.4ex}{$^*$}}

\newtheorem{Pf}{Proof$\!\!$}

\DeclareMathSymbol{\twoheadrightarrow}  {\mathrel}{AMSa}{"10}
\newcommand{\onto}{\protect{~-\!\!\!\twoheadrightarrow}}

\newcounter{letter}
\renewcommand{\theletter}{\rom{(}\alph{letter}\rom{)}}

\newcounter{rnum}
\renewcommand{\thernum}{\rom{(}\roman{rnum}\rom{)}}

\begin{document}

\title[Algebra \& Geometry in the Setting of Regular Algebras]%
{The Interplay of Algebra and Geometry in the \\[2mm]
Setting of Regular Algebras}


\maketitle

\vspace*{0.1in}

\baselineskip15pt

\begin{center}
\renewcommand{\thefootnote}{\fnsymbol{footnote}}
{\sc Michaela Vancliff}\footnote{M.~Vancliff was
supported in part by NSF grants DMS-0900239 and DMS-1302050 and by general %
membership at MSRI.}\\
Department of Mathematics, P.O.~Box 19408,\\
University of Texas at Arlington,\\
Arlington, TX \ 76019-0408\\
{\sf vancliff@uta.edu \quad uta.edu/math/vancliff}
\end{center}

\setcounter{page}{1}
\thispagestyle{empty}

\bigskip
\bigskip

\begin{abstract}
\baselineskip13pt
This article is based on a talk given by the author at MSRI in the
workshop {\it Connections for Women} in January 2013, while being a part
of the program {\it Noncommutative Algebraic Geometry and Representation 
Theory} at MSRI.
One purpose of the exposition is to motivate and describe the geometric 
techniques introduced by M.~Artin, J.~Tate and M.~Van den Bergh in the
1980s at a level accessible to graduate students.  Additionally, some
advances in the subject since the early 1990s are discussed, including a
recent generalization of complete intersection to the noncommutative
setting, and
the notion of graded skew Clifford algebra and its application to classifying 
quadratic regular algebras of global dimension at most three. 
The article concludes
by listing some open problems.
\end{abstract}

\baselineskip16.5pt


\bigskip

\section*{Introduction}

Many non-commutative algebraists in the 1980s were aware of the successful 
marriage of algebra and algebraic geometry in the commutative setting and 
wished to duplicate that relationship in the non-commutative setting.  One such
line of study was the search for a subclass of non-commutative algebras that 
``behave'' enough like polynomial rings that a geometric theory could be 
developed for them. One proposal for such a class of algebras are the {\em 
regular} algebras, introduced in \cite{AS}, that were investigated using new 
geometric techniques in the pivotal papers of M.~Artin, J.~Tate and M.~Van den
Bergh (\cite{ATV1,ATV2}).

About the same time, advances in quantum mechanics in the 20th century 
had produced many new non-commutative algebras on which traditional techniques 
had only yielded limited success, so a need had arisen to find new techniques
to study such algebras (c.f., \cite{FRT,KKO,Sk,Sk2,Sud}).
One such algebra was the Sklyanin algebra, which
had emerged from the study of quantum statistical mechanics (\cite{Sk,Sk2}).
By the early
1990s, T.~Levasseur, S.~P.~Smith, J.~T.~Stafford and others had solved the 
10-year old open problem of completely classifying all the finite-dimensional 
irreducible representations (simple modules) over the Sklyanin algebra, 
and their methods were the geometric techniques developed by Artin, Tate and 
Van den Bergh (\cite{Lev.Smith,Smith.Staff,Smith.Stan}).

Concurrent with the above developments, another approach was considered via 
differential geometry and deformation theory to study the algebras produced by 
quantum physics. That approach is the study of certain non-commutative algebras
via Poisson geometry (c.f., \cite{D}). 
At the heart of both approaches are homological
and categorical techniques, so it is perhaps no surprise that the two 
approaches have much overlap; often, certain geometric objects from one 
approach are in one-to-one
correspondence with various geometric objects from the other approach
(depending on the algebra being studied -- c.f., \cite{V2,V3,V4}).
A survey of recent advances in Poisson geometry may be found 
in~\cite{G}.

Given the above developments, the early 1990s welcomed a new era in the field 
of non-commutative algebra in which geometric techniques took center stage.
Since that time, the subject has spawned many new ideas and directions, as 
demonstrated by the MSRI programs in 2000 and 2013.

\bigskip

This article is based on a talk given by the author in the {\em Connections 
for Women} workshop held at MSRI in January 2013 and it has two objectives.
The first is to motivate and describe the geometric techniques of Artin, Tate 
and Van den Bergh at a level accessible to graduate students, and the second 
is to discuss some developments towards the attempted classification of 
quadratic 
regular algebras of global dimension four, while listing open problems.
An outline of the article is as follows.

Section~1 concerns the motivation and development of the subject, with 
emphasis on quadratic regular algebras of global dimension four. Section~2
discusses constructions of certain types of quadratic regular algebras of
arbitrary finite global dimension, with focus on graded Clifford algebras and 
graded skew Clifford algebras.  This section also discusses a new type of 
symmetry for square matrices called $\mu$-symmetry.  We conclude this section 
by revisiting the classification of quadratic regular algebras of global 
dimension at most three, since almost all such algebras may be formed from 
regular graded skew Clifford algebras.  
In Section~3, we discuss geometric techniques that apply to graded Clifford 
algebras and graded skew Clifford algebras in order to determine
when those algebras are regular. This section also considers the issue of 
complete intersection in the non-commutative setting.  We conclude with 
Section~4 which lists some open problems and related topics.

Although the main objects of study from \cite{ATV1,ATV2} are discussed in
this article, several topics from \cite{ATV1,ATV2} are omitted; for surveys 
of those topics, the reader is referred to \cite{Staff,SVdB} and to 
D.~Rogalski's lecture notes, \cite{R}, from the graduate workshop 
``Noncommutative Algebraic Geometry'' at MSRI in June~2012.

\baselineskip17pt

\begin{ack}
The author is grateful to MSRI for general membership in 2013 and for
providing an office during her visits. Additionally, special thanks are
extended to the organizers of the workshop {\it Connections for Women} in 
January 2013 for their invitation to give a talk in the workshop.
The environment at MSRI was particularly stimulating during the program
{\it Noncommutative Algebraic Geometry and Representation Theory} in Spring 
2013, and for that the author warmly thanks the organizers and participants 
of the program, especially T.~Cassidy, M.~Van den Bergh and J.~T.~Stafford. 
\end{ack}

\bigskip


\section{The Geometric Objects}

In this section, we discuss the motivation and development of the subject, with 
emphasis on quadratic regular algebras of global dimension at most four. 

Throughout this section, $\kk$ denotes an algebraically closed field and,
for any graded algebra~$B$, the span of the homogeneous elements of degree $i$
will be denoted by $B_i$. 

\medskip

\subsection{Motivation}\label{mot}\hfill\\
\indent
Consider the $\kk$-algebra, $S$, on generators
$z_1$, $\ldots$, $z_n$ with defining relations:
$$z_j z_i = \mu_{ij} z_i z_j, \quad \text{ for all distinct } i,~j,$$
where $0 \neq \mu_{ij} \in \kk$ for all $i$, $j$,
and $\mu_{ij} \mu_{ji} = 1$ for all distinct $i$, $j$.
If $\mu_{ij} =1$, for all $i$, $j$, then $S$ is the commutative polynomial 
ring and has a rich subject of algebraic geometry associated with it;
in particular, by the (projective) Nullstellensatz, the points of 
$\PP({S_1}\!\dual)$ 
are in one-to-one correspondence with certain ideals of
$S$ via $(\alpha_1, \ldots , \alpha_n) \leftrightarrow \la \alpha_i z_1 -
\alpha_1 z_i, \ldots , \alpha_i z_n - \alpha_n z_i \ra $, where $\alpha_i
\neq 0$. Before continuing, we first observe that for such an ideal $I$,
the graded module $S/I$ has the property that its Hilbert series is 
$H(t) = 1/(1-t)$ and that $S/I$ is a 1-critical (with respect to GK-dimension) 
graded cyclic module over $S$. 

However, if $\mu_{ij} \neq 1$ for any $i$, $j$, then $S$ still ``feels'' close 
to commutative, and one would expect there to be a way to relate algebraic 
geometry to it.  The geometric objects in \cite{ATV1} are modelled on the 
module $S/I$ above;  instead of using actual points or lines etc, certain 
graded modules are used as follows.

\medskip

\subsection{Points, Lines, etc}\hfill\\[-9mm]
\begin{defn}\label{ptmod}\cite{ATV1} \  
Let $A = \bigoplus_{i= 0}^\infty A_i$ denote an $\N$-graded, connected (meaning 
$A_0 =\kk$) $\kk$-algebra generated by $A_1$ where $\dim(A_1) = n 
< \infty$.  A graded right $A$-module $M = \bigoplus_{i = 0}^\infty M_i$ 
is called a {\it right point module} (respectively, {\it line module}) if:
\begin{enumerate}
\item[(a)] $M$ is cyclic with $M = M_0 A$, and 
\item[(b)] $\dim_{\kk} (M_i) = 1$ for all $i$ (respectively, 
           $\dim_{\kk} (M_i) = i + 1$) for all $i$. 
\end{enumerate}
\end{defn}

\noindent 
If $A$ is the polynomial ring $S$, then the module $S/I$ from \S\ref{mot} is a
point module. In general, one may associate some geometry to point
and line modules as follows. Condition~(a) implies that $A$ maps onto $M$
via $a \mapsto ma$, for all $a \in A$, where $\{m\}$ is a $\kk$-basis for 
$M_0$,
and this map restricts via the grading to a linear map $\theta: A_1 \to M_1$.
Let $K \subset A_1$ denote the kernel of $\theta$. Condition~(b) implies
that $\dim_\kk (K) = n-1$ (respectively, $n-2$), so that $K^\perp \subset
{A_1}\!\dual$ has dimension one (respectively, two). Thus, $\PP(K^\perp)$ is a 
point (respectively, a line) in the geometric space $\PP({A_1}\!\dual)$.

The Hilbert series of a point module is $H(t) = 1/(1-t)$, whereas the
Hilbert series of a line module is $1/(1-t)^2$. Hence, a plane module is
defined as in Definition~\ref{ptmod} but condition~(b) is replaced by the
requirement that the module have Hilbert series $1/(1-t)^3$
(c.f., \cite{ATV1}).  Similarly, one may define $d$-linear modules, where
the definition is modelled on Definition~\ref{ptmod}, but the module has 
Hilbert series $1/(1-t)^{d+1}$ (c.f., \cite{Shelton.V3}).

For many algebras, $d$-linear modules
are $(d+1)$-critical with respect to GK-dimension. This leads to the following 
generalization of a point module.

\begin{defn}\label{bptmod}\cite{CVc} \  
With $A$ as in Definition~\ref{ptmod}, we define a {\em right base-point 
module} over~$A$ to be a graded 1-critical (with respect to GK-dimension)
right $A$-module $M$ such that $M = \bigoplus_{i=0}^\infty M_i = M_0 A$
and $M$ has Hilbert series $H_M(t) = c/(1-t)$ for some $c\in \N$.
\end{defn}

\noindent
If $c =1$ in Definition~\ref{bptmod}, then the module is a point module;
whereas if $c \geq 2$, then the module is called a fat point module (\cite{A}).
The only base-point modules over the polynomial ring are point modules.
On the other hand, in general, the algebra $S$ from \S\ref{mot} can have 
fat point modules,
so fat point modules are viewed as generalizations of points, and this is made 
more precise in \cite{A}.

In \cite{ATV1}, Artin, Tate and Van den Bergh proved that, under certain
conditions, the point modules are parametrized by a scheme; that is, there
is a scheme that represents the functor of point modules. Later, in 
\cite{V.VanR}, this scheme was called the {\it point scheme}.
A decade later, in \cite{Shelton.V3}, it was proved by B.~Shelton and the 
author that (under certain conditions) $d$-linear modules are parametrized by 
a scheme; that is, there is a scheme that represents the functor of $d$-linear 
modules.  If $d = 0$, then this scheme is isomorphic to the point scheme; 
if $d = 1$, the scheme is called the {\it line scheme}.

By factoring out a nonzero graded submodule from a point module, one obtains 
a truncated point module as follows.

\begin{defn}\label{tptmod}\cite{ATV1} \  
With $A$ as in Definition~\ref{ptmod}, we define a {\em truncated right point 
module of length $m$} to be a graded right $A$-module $M = \bigoplus_{i
=0}^{m-1} M_i$ such that $M$ is cyclic, $M = M_0 A$ and $\dim_{\kk} 
(M_i) = 1$ for all $i = 0, \ldots , m-1$.
\end{defn}

\noindent
For many quadratic algebras $A$, there exists a one-to-one correspondence 
between the truncated point modules over~$A$ of length three and the point 
modules over~$A$. Moreover, if the algebra $A$ in Definition~\ref{tptmod} is 
quadratic, then the truncated point modules of length three are in one-to-one 
correspondence with the zero locus in $\PP({A_1}\!\dual)\times
\PP({A_1}\!\dual)$ 
of the defining relations of~$A$. To see this, we fix a $\kk$-basis $\{ x_1,
\ldots , x_n\}$ for $A_1$, and use $T$ to denote the free $\kk$-algebra
on $x_1, \ldots , x_n$, and let $Z \subset \PP({A_1}\!\dual)\times
\PP({A_1}\!\dual)$
denote the zero locus of the defining relations of~$A$. Viewing each
$x_i$ as the $i$'th coordinate function on ${A_1}\!\dual$, let 
$p = (\alpha_i) \in \PP({A_1}\!\dual)$ and $r = (\beta_i) \in
\PP({A_1}\!\dual)$,
where $\alpha_i$, $\beta_i \in \kk$ for all $i = 1, \ldots , n$. Let $M = 
\kk v_0 \oplus \kk v_1 \oplus \kk v_2 $ denote a three-dimensional vector
space that is a $T$-module via the action determined by 
$$
v_0 x_i = \alpha_i v_1, \qquad v_1 x_i = \beta_i v_2, \qquad v_2 x_i = 0,
$$
for all $i$. It follows that $M$ is a truncated point module over $T$ of
length three. If $g \in T_2$, then $v_1 g = 0 = v_2 g$ and 
$v_0 g = g(p, r)v_2$. In particular, if $f\in T_2$ is a defining relation
of $A$, then $Mf = 0$ if and only if $f(p, r) = 0$. Hence, $M$ is an 
$A$-module if and only if $(p, r) \in Z$. This one-to-one correspondence
between $Z$ and truncated point modules of length three
also exists at the level of schemes; the reason being that the scheme~$Z$
represents the functor of truncated point modules of length three. The
method of proof of this is to repeat the preceding argument for a truncated 
point module of length three over $R\otimes_\kk T$ and $R\otimes_\kk A$, where
$R$ is a commutative $\kk$-algebra, together with localization techniques; for
details the reader is referred to \cite[Proposition~3.9]{ATV1}, its proof, and
the paragraph preceding that result. This correspondence will be revisited
in \S\ref{gldim4}.

\bigskip

For completeness, we finish this subsection with some technical definitions 
that play minor roles throughout the text.  The reader is referred to
\cite{Lev,Lev.Smith} for details and for results concerning algebras 
satisfying these definitions.

\begin{defn} \cite[Definition~2.1]{Lev}
A noetherian ring $B$ is called {\em Auslander-regular} (respectively, 
{\em Auslander-Gorenstein}) if 
\begin{enumerate}
\item[(a)] the global homological dimension (respectively, (left and right) 
           injective dimension) of~$B$ is finite, and
\item[(b)] every finitely generated $B$-module~$M$ satisfies the {\em 
           Auslander condition}, namely, for every $i \geq 0$ and for 
	   every $B$-submodule~$N$ of Ext$^i_B(M,\, B)$, we have $j(N) \geq 
	   i$, where
           $j(N) = \text{inf} \{ \ell : \text{Ext}^\ell_B(N,\, B) 
	    \neq 0\}.$
\end{enumerate}
\end{defn}

\begin{defn} \cite[Definition~5.8]{Lev}
A noetherian $\kk$-algebra $B$ of integral GK-dimension~$n$ 
satisfies the {\em Cohen-Macaulay property} if \ GKdim$(M) + j(M) = n$ \  
for all nonzero finitely generated $B$-modules~$M$.
\end{defn}

\bigskip

\subsection{Regular Algebras}\hfill\\
\indent
The goal of \cite{ATV1} was to classify, in a user-friendly way, the generic 
regular algebras of global dimension three that were first analysed in 
\cite{AS}. In \cite{ATV1}, such algebras were shown to be noetherian by using 
the geometric techniques developed in \cite{ATV1}. Regular algebras are often
viewed as non-commutative analogues of polynomial rings and are defined as 
follows.

\begin{defn}\cite{AS}\label{regalg} \  
A finitely generated, $\N$-graded, connected $\kk$-algebra 
$A = \bigoplus_{i= 0}^\infty A_i$, generated by $A_1$, is {\em regular}
(or {\em AS-regular}) of global dimension~$r$ if
\begin{enumerate}
\item[(a)]  it has global homological dimension~$r < \infty$, and 
\item[(b)]  it has polynomial growth (i.e., there exist positive real 
            numbers $c$ and $\delta$ such that $\dim_\kk (A_i) \leq
	    c i^\delta$ for all $i$), and  
\item[(c)]  it satisfies the {\em Gorenstein condition}, namely, a 
            minimal projective resolution of the left trivial module ${}_A\kk$
	    consists of finitely generated modules and dualizing this 
	    resolution yields a minimal projective resolution of the right 
	    trivial module $\kk_A[e]$, shifted by some degree $e$.
\end{enumerate}
\end{defn}

\noindent
Although all three conditions in Definition~\ref{regalg} are satisfied by
the polynomial ring, the main reason a regular algebra is viewed as a
non-commutative analogue of a polynomial ring is due to condition~(c), since
it imposes a symmetry condition on the algebra that replaces the symmetry
condition of commutativity. The reader should note that, in the literature, 
(c)~is sometimes replaced by an
equivalent condition that makes the symmetry property less obvious; namely,
Ext$^i_A({}_A\kk, A) \cong \delta^i_r\, \kk_A[e]$, where $\delta^i_r$ is the 
Kr\"onecker-delta symbol.
An $\N$-graded connected $\kk$-algebra that is generated by degree-1 elements
and which is Auslander-regular with polynomial growth is AS-regular 
(\cite{Lev}).
For a notion of regular algebra where the algebra is not generated 
by degree-1 elements, see \cite{C1,C2,Steph2,Steph1,Steph3,Steph4}.

\begin{examples}\label{eg1}\hfill\\
\indent (a) The algebra $S$ from \S\ref{mot} is regular.\\
\indent (b) If $\kk =\C$, then many algebras from physics are regular.
           In particular, homogenizations of universal enveloping algebras
	   of finite-dimensional Lie algebras, the coordinate ring of quantum 
	   affine $n$-space,
	   the coordinate ring of quantum $m \times n$ matrices, and the
	   coordinate ring of quantum symplectic $n$-space are all regular
	   (\cite{Bruyn.Smith,Bruyn.VdB,Lev.Staff})\\
\indent (c) 
If the global dimension of a regular algebra is one, then the algebra is
the polynomial ring on one variable.  However, by \cite{AS}, if the global
dimension is two, then there are two types of such algebra as follows. For
both types, the
algebra has two generators, $x$, $y$, of degree one and one defining 
relation~$f$, where either $f = xy - yx - x^2$ ({\em Jordan plane}) or 
$f = xy - q yx$ ({\em quantum affine plane}), where $q\in \kk$ can be any 
nonzero scalar.
\end{examples}

However, if the global dimension is three, then the situation is much richer;
some of the algebras are quadratic with three generators and three defining 
relations, whereas the rest have two generators and two cubic relations
(\cite{AS}).
Such algebras that are generic are classified in \cite{ATV1} according to their
point schemes, and in all cases, the point scheme is the graph of an
automorphism $\sigma$. Moreover, the algebra is a finite module over its center
if and only if $\sigma$ has finite order.

\medskip

\subsection{Global Dimension Four} \label{gldim4} \hfill\\
\indent
Although many regular algebras of global dimension four have been extensively 
studied, there is no classification yet. 
Recently, the progress towards classifying non-quadratic regular algebras of 
global dimension four made good headway via the work in \cite{LPWZ,RZ}. 
However, quadratic regular algebras of global dimension four constitute most 
of the regular algebras of global dimension four, so their attempted 
classification is one of the motivating problems that drives the subject 
forward. 
We end this section by summarizing some key results for this latter case; 
in this setting, the algebra has four generators and six relations.

In unpublished work, Van den Bergh proved in the mid-1990s that any quadratic
(not necessarily regular) algebra $A$ on four generators with six generic 
defining relations has twenty (counted with multiplicity) nonisomorphic 
truncated point modules of length three.  Hence, $A$ has at most twenty 
nonisomorphic point modules.  He also proved that if, additionally, $A$ is
Auslander-regular of global dimension four, then $A$ has a 1-parameter family 
of line modules. For lack of a suitable reference, we outline the proof of 
these results.  Let $M(4,\, \kk)$ denote the space of $4 \times 4$ matrices 
with entries in $\kk$. For the first result, we write points of
$\PP({A_1}\!\dual)$
as columns and, by mapping  $(a, b) \in \PP({A_1}\!\dual) \times
\PP({A_1}\!\dual)$
to the matrix $ab^T \in M(4,\, \kk)$, we have that $\PP({A_1}\!\dual) \times 
\PP({A_1}\!\dual)$ is isomorphic to the scheme $\Omega_1$ of rank-1 elements in 
$\PP(M(4,\, \kk))$. Correspondingly, the defining relations of $A$ map to 
homogeneous degree-1 polynomial functions on $M(4,\, \kk)$, and their zero 
locus $Z' \subset \PP(M(4,\, \kk))$ can be identified with a $\PP^9$. With 
these identifications, the zero locus $Z \subset \PP({A_1}\!\dual) \times 
\PP({A_1}\!\dual)$ of the defining relations of $A$ is isomorphic to 
$\Omega_1 \cap Z' \subset \PP(M(4,\, \kk))$. Since $\Omega_1$ has dimension 
six and degree twenty, $\dim(Z) \geq 6+9-15 = 0$, and, by B\'ezout's Theorem, 
$\deg(Z)=20$. Hence, generically, $Z$ is finite with twenty points, so the 
first result follows by using the discussion after Definition~\ref{tptmod}.
For the second result, we identify $A_1 \otimes_\kk A_1$ with $M(4,\, \kk)$, 
and the assumption on regularity allows
the application of \cite[Proposition~2.8]{Lev.Smith}, so that the line modules
are in one-to-one correspondence with the elements in the span of the defining
relations of $A$ that have rank at most two. In particular, we compute
$\dim(\Omega_2 \cap \Delta)$ in $\PP(M(4,\, \kk))$, where $\Omega_2$ denotes
the elements in $\PP(M(4,\, \kk))$ of rank at most two and $\Delta$ denotes 
the projectivization of the image in $\PP(M(4,\, \kk))$ of the span of the 
defining relations of $A$. Since $\Delta \cong \PP^5$ and $\dim(\Omega_2) =11$,
the dimension is thus at least equal to $11 + 5 - 15 = 1$, so, generically, 
$A$ has a 1-parameter family of line modules.

In spite of Van den Bergh's work, it was still not clear that a regular 
algebra
satisfying the hypotheses from the preceding paragraph could have both a 
finite point scheme (especially one of cardinality twenty) and a 1-dimensional
line scheme simultaneously.  However, in \cite{VVW}, the author proved with 
Van Rompay and Willaert, in the mid-1990s,
that there exists a quadratic regular algebra of global dimension four on four 
generators with six defining relations that has exactly one point module (up 
to isomorphism) and a 1-parameter family of line mods.

Some years later, in 2000, Shelton and the author proved in \cite{Shelton.V3} 
that if a quadratic algebra on four generators with six defining relations
has a {\em finite} scheme of truncated point modules of length three, then that
scheme determines the defining relations of the algebra. One should note that 
this result assumes no hypothesis of regularity nor of any other homological 
data.  Moreover, by \cite{VanR}, this result is false in general if the scheme
is infinite, even if the algebra is assumed to be regular and noetherian.

Shelton and the author also proved in \cite{Shelton.V3} that if a quadratic 
regular algebra of global dimension four (that satisfies a few other 
homological conditions) has four generators and six defining relations and a 
1-dimensional line scheme, then that scheme determines the defining relations 
of the algebra.

These last two results are counter-intuitive, since they seem to be saying
that if the point scheme (respectively, line scheme) is {\em as small as 
possible}, then the defining relations can be recovered from it.

However, by the start of 2001, it was still unclear whether or not any 
quadratic regular 
algebra exists that has global dimension four, four generators, six defining 
relations, exactly twenty nonisomorphic point modules and a 1-dimensional line 
scheme. Fortunately, this was resolved by Shelton and Tingey in 
\cite{S.Tingey} in 2001 in the affirmative. Sadly, their method to
produce their example used much trial and error on a computer, which they 
and others were unable to duplicate to produce more examples. This hurdle
likely had a negative impact on the development of the subject, since it
is difficult to make conjectures if there is only one known example.
Hence, a quest began to find an algorithm to construct such algebras, but
it was another several years before this situation was remedied, and that
is discussed in the next section.

\bigskip
\bigskip


\section{Graded Clifford Algebras, Graded Skew Clifford\\[0mm]
Algebras and Quantum Planes}\label{sec2}

This section describes a construction of a certain type of regular algebra of 
arbitrary finite global dimension; such an algebra is called a graded skew 
Clifford algebra as it is modelled on the construction of a graded Clifford 
algebra.  If the global dimension is four, then this construction is 
able to produce regular algebras that have the desired properties described at 
the end of the previous section. We conclude this section by revisiting the 
classification of quadratic regular algebras of global dimension three, and
show that almost all such algebras may be obtained from regular graded skew 
Clifford algebras. 

We continue to assume that $\kk$ is algebraically closed; we additionally 
assume char$(\kk) \neq 2$. We write $M(n,\,\kk)$ for the space of 
$n \times n$ matrices with entries in $\kk$, and $M_{ij}$ for the entry
in the $n \times n$ matrix~$M$ that is in row~$i$ and column~$j$.

\medskip

\subsection{Graded Clifford Algebras}
\begin{defn}\label{gca}\label{gcas}\cite{Aubry.Lemaire,L}
Let $M_1 , \ldots , M_n \in M(n,\,\kk)$ denote symmetric matrices.  A {\em 
graded Clifford algebra} (GCA) is the $\kk$-algebra~$C$ on degree-one 
generators 
$x_1, \ldots , x_n$ and on degree-two generators $y_1, \ldots , y_n$ with 
defining relations given by:
\begin{enumerate}
\item[(i)] (degree-2 relations) \  
            $\ds x_i x_j + x_j x_i = \sum_{k=1}^n (M_k)_{ij} \, y_k$ for all 
	    $i, j = 1, \ldots , n$, and\\[-3mm]
\item[(ii)] degree-3 and degree-4 relations that guarantee 
           $y_k$ is central in $C$ for all $k = 1, \ldots , n$. 
\end{enumerate}
\end{defn}
\noindent
In general, GCAs need not be quadratic nor regular, as demonstrated by
the next example.
\begin{example}\label{gcaeg}
\baselineskip21pt
Let $\scriptsize M_1 = \begin{bmatrix} 2 & -1 \\ -1 & 0  \end{bmatrix}$ and 
$\scriptsize M_2 = \begin{bmatrix} 0 & -1 \\ -1 & 2  \end{bmatrix}$. 
The corresponding GCA is the $\kk$-algebra on degree-one generators
$x_1,\ x_2$ with defining relations 
$$
x_1 x_2  + x_2 x_1  = -x_1 ^2 - x_2^2, \qquad x_1 ^2 x_2 = x_2 x_1 ^2,
$$
so this algebra is not quadratic nor regular (as $(x_1+x_2)^2 = 0$). For more
details on this algebra, the reader may consult \cite[Example~2.4]{V5}.
\end{example}
\baselineskip17pt

GCAs $C$ are noetherian by \cite[Lemma~8.2]{ATV1}, since 
$\dim_\kk(C/\la y_1 , \ldots , y_n \ra) < \infty$. Moreover, 
since each matrix $M_k$ in the definition is symmetric, we may associate a 
quadratic form to~$M_k$, and thereby associate a quadric in $\PP^{n-1}$
to $M_k$ for each $k$. This means that for each GCA, as in 
Definition~\ref{gca}, there is an associated quadric system $\Q$ in 
$\PP^{n-1}$.  Quadric systems are said to be base-point free if they yield a
complete intersection; that is, the intersection of all the quadrics in the 
quadric system is empty.  Although Example~\ref{gcaeg} demonstrates that a 
GCA need not be quadratic nor regular, if $\Q$ is base-point free,
it determines these properties of the associated GCA as follows.

\begin{thm}\label{gcathm}\cite{Aubry.Lemaire,L} \ 
The GCA $C$ is quadratic, Auslander-regular of global dimension~$n$ and 
satisfies the Cohen-Macaulay property with Hilbert series $1/(1-t)^n$ if and 
only if the associated quadric system is base-point free; in this case, 
$C$~is regular and a domain.
\end{thm}

In spite of this result, regular GCAs of global dimension four are not 
candidates for generic quadratic regular algebras of global dimension four, 
since, although their point schemes can be finite (\cite{Steph.V2,VVW}), the 
symmetry of their relations prevents their line schemes from having dimension 
one (\cite{Shelton.V3}).  The standard argument to prove this for a
quadratic regular GCA~$C$ of global dimension four exploits the symmetry of 
the defining relations of $C$ to move the computation of \S\ref{gldim4}
inside $\PP(W)$, where $W$ is the 10-dimensional subspace of $M(4,\, \kk)$ 
consisting of all symmetric matrices. Hence, using the notation from 
\S\ref{gldim4}, $\Delta \subset \PP(W)$ and 
the line modules are parametrized by $(\Omega_2\cap
\PP(W)) \cap \Delta \subset \PP(W)$; thus the dimension is at least 
$6 + 5 - 9$, so it is at least two.

Hence, a modification of the definition of GCA is desired in such a way
that enough symmetry is retained so as to allow an analogue of 
Theorem~\ref{gcathm} to hold, while, at the same time, losing some symmetry 
so that the line scheme might have dimension one.

\bigskip

\subsection{Graded Skew Clifford Algebras}\label{gscas}\hfill

\indent
In order to generalize the notion of GCA and to have a result analogous to
Theorem~\ref{gcathm}, we need to generalize the notions of symmetric matrix
and quadric system and make use of normalizing sequences.
For any $\N$-graded $\kk$-algebra $B$, a sequence $\{g_1, \ldots , g_m\}$
of homogeneous elements of positive degree is called {\em normalizing} if 
$g_1$ is a normal element in $B$ and, for each $k = 1, \ldots , m-1$, the image 
of $g_{k+1}$ in $B/\la g_1, \ldots, g_k\ra$ is a normal element. 

We write $\kk^\times$ for $\kk\setminus\{0\}$.

\begin{defn}\label{gsca}\cite{CV}\\
\indent
(a) Let $\mu \in M(n,\,\kk^\times)$ satisfy $\mu_{ij} \mu_{ji} = 1$ for all
distinct $i$, $j$. We say a matrix $M \in M(n,\,\kk)$ is $\mu$-{\em symmetric} 
if $M_{ij} = \mu_{ij}M_{ji}$ for all $i, j = 1, \ldots , n$. We write
$M^\mu (n,\,\kk)$ for the subspace of $M(n,\,\kk)$ consisting of all 
$\mu$-symmetric matrices.\\
\indent
(b) Fix $\mu$ as in (a) and additionally assume $\mu_{ii} =1$ for all $i$.
Let $M_1, \ldots , M_n \in M^\mu (n,\,\kk)$.
A {\em graded skew Clifford algebra} (GSCA) associated to $\mu$ and $M_1$, 
$\ldots ,$ $M_n$ is a graded $\kk$-algebra
$A = A(\mu, M_1, \ldots , M_n)$
on degree-one generators $x_1, \ldots , x_n$ and on degree-two generators 
$y_1, \ldots , y_n$ with defining relations given by:
\begin{enumerate}
\item[(i)] (degree-2 relations)
           $\ds x_i x_j + \mu_{ij} x_j x_i = \sum_{k=1}^n (M_k)_{ij} y_k$
           for all $i, j = 1, \ldots , n$, and\\[-3mm]
\item[(ii)] degree-3 and degree-4 relations that guarantee the existence of a 
           normalizing sequence $\{ y_1', \ldots , y_n'\}$ that spans 
	   $\sum_{k=1}^n\kk y_k$. 
\end{enumerate}
\end{defn}

\noindent
Clearly, symmetric matrices and skew-symmetric matrices are $\mu$-symmetric 
matrices for appropriate $\mu$, and GCAs are GSCAs. Moreover, by 
\cite[Lemma~8.2]{ATV1}, GSCAs $A$ are noetherian since 
$\dim_\kk(A/\la y_1 , \ldots , y_n \ra) < \infty$.  
Furthermore, in Definition~\ref{gsca}(b)(i), for all $i$, $j$, the 
$ji$-relation can be deduced from the $ij$-relation by the $\mu$-symmetry of 
the $M_k$.

\begin{examples}\label{gscaeg1}\hfill\\
\indent (a)  With $\mu$ as in Definition~\ref{gsca}(b),
skew polynomial rings on generators \ $x_1, \ldots , x_n$ \ with
relations \ $x_i x_j = - \mu_{ij}x_j x_i$, \ for all $i \neq j$, are GSCAs.
\\[3mm]
\indent (b)  (Quantum Affine Plane) \  
Let $n =2$, and $M_1 = \scriptsize \begin{bmatrix} 2 & 0 \\ 0 & 0
\end{bmatrix}$ and 
$M_2 = \scriptsize \begin{bmatrix} 0 & 0 \\ 0 & 2 \end{bmatrix}$.
The degree-2 relations of $A(\mu,\,M_1,\,M_2)$ have the form:
\[
2 x_1^2 = 2 y_1, \qquad 2 x_2^2 = 2 y_2, \qquad 
 x_1 x_2 + \mu_{12} x_2 x_1 = 0,
\]
so that \  
$\kk \la x_1, x_2 \ra / \la x_1 x_2 + \mu_{12} x_2 x_1 \ra \ \onto 
\ A(\mu,\,M_1,\,M_2)$. By Theorem~\ref{gscathm} below, this map
is an isomorphism (see Examples~\ref{gscaeg2}(a)).\\[3mm]
\indent (c)  (``Jordan'' Plane) \  
Let $n =2$, and $M_1 = \scriptsize\begin{bmatrix} 2 & 1 \\ \mu_{21} & 0
\end{bmatrix}$ 
and $M_2 = \scriptsize\begin{bmatrix} 0 & 0 \\ 0 & 2 \end{bmatrix}$.
The degree-2 relations of $A(\mu,\,M_1,\,M_2)$ have the form:
\[
2 x_1^2 = 2 y_1, \qquad 2 x_2^2 = 2 y_2, \qquad
 x_1 x_2 + \mu_{12} x_2 x_1 = y_1 = x_1^2,
\]
so that \  
$\kk \la x_1, x_2 \ra / \la x_1 x_2 + \mu_{12} x_2 x_1 - x_1^2\ra \ \onto 
\ A(\mu,\,M_1,\,M_2)$. By Theorem~\ref{gscathm} below, this map
is an isomorphism (see Examples~\ref{gscaeg2}(b)). Depending on the choice
of $\mu_{12}$, this family of examples contains the Jordan plane and some 
quantum affine planes.\\[3mm]
\indent (d) The quadratic regular algebra of global dimension four found by 
Shelton and Tingey in 2001, in \cite{S.Tingey}, and discussed above in
\S\ref{gldim4}, that has exactly twenty nonisomorphic point modules and a 
1-dimensional line scheme is a GSCA (\cite{CV}).
\end{examples}

One can associate a non-commutative ``quadric'' to each $\mu$-symmetric matrix
$M_k$ and, in so doing, there is 
also a notion of ``base-point free''.  These ideas are discussed in
\S\ref{ncci} below,
and yield a generalization of Theorem~\ref{gcathm} as follows.

\begin{thm}\label{gscathm}\cite{CV} \  
The GSCA $A$ is quadratic, Auslander-regular of global dimension~$n$ and 
satisfies the Cohen-Macaulay property with Hilbert series $1/(1-t)^n$ if and 
only if the associated quadric system is normalizing and base-point free; in 
this case, $A$~is regular and a domain and uniquely determined, up to 
isomorphism, by the data $\mu$, $M_1, \ldots , M_n$.
\end{thm}

Theorem~\ref{gscathm} allowed the production in \cite{CV} of many algebras
that are candidates for generic quadratic regular algebras of global dimension 
four. In particular, there exist quadratic regular GSCAs of global dimension 
four on four generators with six defining relations that have exactly twenty 
nonisomorphic point modules and a 1-dimensional line scheme.

It is an open problem to describe the 1-dimensional line schemes of the 
regular GSCAs of global dimension four in \cite{CV} that have exactly twenty 
nonisomorphic point modules.

By Examples~\ref{eg1}(c) and \ref{gscaeg1}(b)(c), the regular algebras of
global dimension at most two are GSCAs, and, by \S\ref{qqp}, almost all
quadratic regular algebras of global dimension three are determined 
by GSCAs, so GSCAs promise to be very helpful in the classification of all 
quadratic regular algebras of global dimension four.

\bigskip

\subsection{Quadratic Quantum Planes}\label{qqp}\hfill

\indent
In the language of \cite{A}, a regular algebra of global dimension three
that is generated by degree-1 elements is sometimes called a {\it quantum 
plane} or {\it quantum projective plane} or a {\em quantum} $\PP^2$.  
The classification of the {\it generic} quantum planes is in
\cite{AS,ATV1,ATV2}. In this subsection, we summarize the results of
\cite{NVZ}, in which all {\em quadratic} quantum planes are classified by
using GSCAs.

We continue to assume that $\kk$ is algebraically closed, but its
characteristic is arbitrary unless specifically stated otherwise.

Let $D$ denote a quadratic quantum plane and let $X \subset \PP^2$ denote its 
point scheme.  By \cite[Proposition~4.3]{ATV1} and \cite[Lemma~2.1]{NVZ}, 
there are, in total, four cases to consider:
\begin{itemize}
\item $X$ contains a line, or 
\item $X$ is a nodal cubic curve in $\PP^2$, or 
\item $X$ is a cuspidal cubic curve in $\PP^2$, or 
\item $X$ is a (nonsingular) elliptic curve in $\PP^2$.
\end{itemize}

\begin{thm}\cite{NVZ} \  
Suppose char$(\kk) \neq 2$. If $X$ contains a line, then either $D$ is a
twist, by an automorphism, of a GSCA, or $D$ is a twist, by a twisting system, 
of an Ore extension of a regular GSCA of global dimension two.
\end{thm}

\begin{thm}\cite{NVZ} \  
If $X$ is a nodal cubic curve, then \ $D$ is isomorphic to a $\kk$-algebra
on generators $x_1$, $x_2$, $x_3$ with defining relations:
\[
\lambda x_1 x_2 = x_2 x_1, \qquad
\lambda x_2 x_3 = x_3 x_2 - x_1^2, \qquad
\lambda x_3 x_1 = x_1 x_3 - x_2^2,
\tag{$*$}
\]
where $\lambda \in \kk$ and $\lambda^3 \notin \{0,\, 1\}$.
Conversely, for any such $\lambda$, any quadratic algebra with defining 
relations \thetag{$*$} is a quantum plane and its point scheme is a nodal cubic 
curve in $\PP^2$. Moreover, if char$(\kk) \neq 2$, then
$D$ is an Ore extension of a regular GSCA of global dimension two;
in particular, if $\lambda^3 = -1$, then $D$ is a GSCA. 
\end{thm}

\begin{thm}\cite{NVZ} \  
If\/ $\char(\kk) = 3$, then $X$ is not a cuspidal cubic curve in $\PP^2$.
If\/ $\char(\kk) \neq 3$ and if $X$ is a cuspidal cubic curve in $\PP^2$, 
then \ $D$ is isomorphic to a $\kk$-algebra on generators $x_1$, $x_2$, $x_3$ 
with defining relations:
\[
x_1 x_2 = x_2 x_1 + x_1^2, \quad
x_3 x_1 = x_1 x_3 + x_1^2 + 3 x_2^2, \quad
x_3 x_2 = x_2 x_3 - 3 x_2^2 - 2 x_1 x_3 - 2 x_1 x_2.
\tag{$\dag$}
\]
Moreover, any quadratic algebra with defining relations \thetag{$\dag$} is a
quantum plane; it has point scheme given by a cuspidal cubic curve in
$\PP^2$ if and only if\/ $\char(\kk) \neq 3$.
If\/ $\char(\kk) \neq 2$, then any quadratic algebra with defining relations 
given by \thetag{$\dag$} is an Ore extension of a regular GSCA of global
dimension two.
\end{thm}

It remains to discuss the case that $X$ is an elliptic curve.
In \cite{AS,ATV1}, such algebras are classified into types A, B, E, H,
where some members of each type might not have an elliptic curve as their
point scheme, but a generic member does.

\begin{thm}\cite{NVZ} \   
Suppose that char$(\kk) \neq 2$ and that $X$ is an elliptic curve.
\begin{enumerate}
\item[(a)] Quadratic quantum planes of type H are GSCAs. 
\item[(b)] Quadratic quantum planes of type B are GSCAs. 
\item[(c)] As in \cite{AS,ATV1}, a quadratic quantum plane $D$ of type A is 
           given by a $\kk$-algebra on generators $x$, $y$, $z$ with defining 
	   relations: 
	   \[ a xy + b yx + c z^2 = 0, \quad 
	   a yz + b zy + c x^2 = 0, \quad 
	   a zx + b xz + c y^2 = 0,\]
	   where $a$, $b$, $c \in \kk^\times$, 
	   $(3abc)^3 \neq (a^3 + b^3 + c^3)^3$,  char$(\kk) \neq 3$,
	   and either $a^3 \neq b^3$,  or $a^3 \neq c^3$, or 
	   $b^3 \neq c^3$. In the case that $a^3 = b^3 \neq c^3$, $D$ is
	   a GSCA; whereas in the case $a^3 \neq b^3 = c^3$ $($respectively, 
	   $a^3 = c^3 \neq b^3$$)$, $D$ is a twist, by an automorphism, of a 
	   GSCA.
\end{enumerate}
\end{thm}

In (c) of the last result, the case that $a^3 \neq b^3\neq c^3 \neq a^3$ is 
still open. Moreover, the case when $D$ is of type~E is still open, but this
case only consists of one algebra, up to isomorphism and anti-isomorphism.
However, both type A and type E have the property that the Koszul dual of~$D$
is a quotient of a regular GSCA; so, in this sense, such algebras are
weakly related to GSCAs.

\bigskip


\section{Complete Intersections}\label{secci}

In this section, we define the geometric terms used in Theorem~\ref{gscathm}.
That discussion leads naturally into a consideration of a notion of 
non-commutative complete intersection that mimics the commutative definition. 

We continue to assume that the field $\kk$ is algebraically closed.

\bigskip

\subsection{Commutative Complete Intersection and Quadric Systems}\hfill

\indent
Let $R$ denote the commutative polynomial ring on $n$ generators of degree one.
If $f_1, \ldots , f_m$ are homogeneous elements of $R$ of positive degree,
then $\{ f_1, \ldots , f_m\}$ is a regular sequence in $R$ if and only if 
GKdim$(R/\la f_1, \ldots , f_k \ra) = n - k \geq 0$, for all $k = 1, . . . ,m$.
Geometrically, this corresponds to the zero locus in $\PP({R_1}\!\dual)$ of the 
ideal $J_k = \la f_1, \ldots , f_k \ra$ having dimension $n-1-k \geq -1$ for 
all $k$.  If $\{ f_1, \ldots , f_m\}$ is a regular sequence, then the zero 
locus of $J_m$ (respectively, $R/J_m$) is called a {\em complete 
intersection} (c.f., \cite{E}).

In the setting of \S\ref{gcas}, a quadric system $\Q$ is associated to 
symmetric matrices $M_1, \ldots , M_n$. In that setting, $\Q$ corresponds to a 
regular sequence in $R$ if and only if $\Q$ is a complete intersection, that 
is, if and only if $\Q$ has no base points (a base point is a point that lies 
on all the quadrics in $\Q$).  A non-commutative analogue of this is needed
for Theorem~\ref{gscathm}. 

\bigskip

\subsection{Non-Commutative Complete Intersection and Quadric Systems}%
\label{ncci} \hfill

\indent
The following result uses the notion of base-point module defined in 
Definition~\ref{bptmod}.

\begin{prop}\cite{CV,CVc}\label{ci1}
Let $S$ denote the skew polynomial ring from \S\ref{mot}, and let $f_1,
\ldots , f_n$ denote homogeneous elements of $S$ of positive degree.
If $\{f_1, \ldots , f_n\}$ is a normalizing sequence in $S$, then the following
are equivalent:
\begin{enumerate}
\item[(a)] $\{ f_1, \ldots , f_n\}$ is a regular sequence in $S$,
\item[(b)] $\dim_\kk(S/\la f_1, \ldots , f_n \ra) < \infty$,
\item[(c)] for each $k = 1, \ldots , n$, we have 
           GKdim$(S/\la f_1, \ldots , f_k \ra) = n-k$,
\item[(d)] the factor ring $S/\la f_1, \ldots , f_n \ra$ has no right 
           base-point modules,
\item[(e)] the factor ring $S/\la f_1, \ldots , f_n \ra$ has no left 
           base-point modules.
\end{enumerate}
\end{prop}
\noindent
Such a sequence $\{ f_1, \ldots , f_n\}$ (respectively, $S/\la f_1, \ldots
, f_n \ra$) satisfying the equivalent conditions (a)-(e) from
Proposition~\ref{ci1} is called a {\em complete intersection} in~\cite{CVc}.

In the setting of \S\ref{gscas}, one associates $S$ to the GSCA by
using $\mu$. The isomorphism $M^\mu (n,\,\kk) \to S_2$ defined by 
$M \mapsto (z_1, \ldots , z_n) M (z_1, \ldots , z_n)^T$ associates  
a quadric system $\Q$ to the $\mu$-symmetric matrices $M_1, \ldots , M_n$; 
that is, $\Q$ is the span in $S_2$ of the images of the $M_k$ under this 
map. If $\Q$ is given by a normalizing sequence in $S$, then it is called
a {\em normalizing quadric system}.
By Proposition~\ref{ci1}, if $\Q$ is normalizing, then it 
corresponds to a regular sequence in $S$ if and only if 
it is a complete intersection, that is, if and only if $S/\la \Q \ra$
has no right (respectively, left) base-point modules; this is the meaning
of {\em base-point free} in Theorem~\ref{gscathm}.

\begin{examples}\label{gscaeg2}\hfill\\
\indent (a) \cite{CV} \   
We revisit the quantum affine plane from Examples~\ref{gscaeg1}(b),
where $n = 2$. In that case, $M_i \mapsto q_i = 2 z_i^2 \in S_2$, for $i = 
1,\, 2$. The sequence $\{q_1,\, q_2\}$ is normalizing in $S$ and $\dim(S/\la
q_1,\, q_2 \ra) < \infty$. Thus, by Proposition~\ref{ci1}, the corresponding
quadric system is base-point free.\\[1mm]
\indent (b) \cite{CV} \   
For Examples~\ref{gscaeg1}(c), $n =2$  and $M_1 \mapsto 
q_1 = 2 (z_1^2 + z_1 z_2)$  and $M_2 \mapsto q_2 = 2 z_2^2$. Here, the
sequence $\{ q_2,\, q_1\}$ is normalizing in $S$  and 
$\dim(S/\la q_2,\, q_1 \ra) < \infty$, so by Proposition~\ref{ci1}, the 
corresponding quadric system is base-point free.
\end{examples}

Proposition~\ref{ci1} has recently been extended in \cite{V5} to a family of 
algebras that contains the skew polynomial ring $S$ from \S\ref{mot}. 
In particular, an analogue of Proposition~\ref{ci1} holds
for regular GSCAs, many quantum groups, and homogenizations of
finite-dimensional Lie algebras.

\begin{thm}\cite{V5} \label{ci} \  
Let $A= \bigoplus_{i=0}^\infty A_i$ denote a connected, $\N$-graded
$\kk$-algebra that is generated by $A_1$.
Suppose $A$ is Auslander-Gorenstein of finite injective dimension
and satisfies the Cohen-Macaulay property, and that there exists a
normalizing
sequence $\{ y_1, \ldots , y_\nu\}\subset A\setminus \kk$ consisting
of homogeneous elements such that
GKdim$(A/\la y_1, \ldots , y_\nu \ra ) = 1$. If GKdim$(A) = n \in \N$, and 
if $F = \{f_1, \ldots , f_n\}\subset A\setminus \kk^\times$ is a normalizing 
sequence of homogeneous elements, then the following are equivalent:
\begin{enumerate}
\item[(a)] $F$ is a regular sequence in $A$,
\item[(b)] $\ds \dim_\kk(A/\la F\ra) < \infty$,
\item[(c)] for each $k = 1, \ldots , n$, we have
           GKdim$\left( A/\la f_1 , \ldots , f_k\ra \right) = n-k$, 
\item[(d)] the factor ring $A/\la F\ra$ has no right base-point modules,
\item[(e)] the factor ring $A/\la F\ra$ has no left base-point modules.
\end{enumerate}
\end{thm}

The reader should note that other notions of complete intersection abound in 
the literature, with most emphasizing a homological approach, such as the 
recent work in~\cite{KKZ}.

\bigskip


\section{Conclusion}

In this section, we list some open problems and related topics.
The open problems are not listed in any particular order in
regards to difficulty, and many challenge levels are included, with some
quite computational in nature, and so accessible to junior researchers. 

\subsection{Some Open Problems}\hfill\\
\indent
1.
As stated at the end of \S\ref{sec2}, it is still open whether or not quadratic
quantum planes of type~A with  $a^3 \neq b^3\neq c^3 \neq a^3$ are 
directly related to GSCAs; the analogous problem is also open for type~E.

2.
Is it possible to classify cubic quantum planes by using GSCAs, or by using 
an appropriate analogue of a GSCA?

3.
Is it possible to classify quadratic regular algebras of global dimension
four by using GSCAs? Presumably, such a classification will use both the
point scheme and the line scheme.

4.
Can standard results on commutative quadratic forms and quadrics be extended 
to non-commutative quadratic forms and quadrics? For example, P.~Veerapen
and the author have extended, in \cite{V.PV1}, the notion of rank of a 
(commutative) quadratic 
form to non-commutative quadratic forms on $n$~generators, where $n = 2,\ 3$;
can this be done for $n \geq 4$? 

5.
Can results concerning GCAs be carried over to GSCAs? In particular, Veerapen 
and the author applied their aforementioned generalization of rank to GSCAs in 
a way that is analogous to that used for the traditional notion of rank with 
GCAs in \cite{VVW}. They proved in \cite{V.PV2} that various results in 
\cite{VVW} concerning point modules over GCAs apply to point modules over 
GSCAs. 

6.
Can standard results concerning symmetric matrices be extended or generalized 
to $\mu$-symmetric matrices?

7. 
Can the results in \cite{V5}, mentioned above at the end of \S\ref{secci}, on 
complete intersections be extended to an even larger family of algebras than 
those considered in \cite{V5}?

8.
By combining results in \cite{CV} and \cite{Steph.V2}, it is known that
regular GSCAs of global dimension four can have exactly $N$~nonisomorphic
point modules, where $N \notin \{2,\,19\}$; it is not yet known if 
$N \in \{2,\,19\}$ is possible. In fact, by \cite{Steph.V1}, $N = 2$ is 
possible if the algebra is quadratic and regular of global dimension four but 
is not a GSCA, but it is not known if $N = 19$ is possible, even if the 
algebra is not a GSCA.

9. 
What is the line scheme of some known quadratic regular algebras of global
dimension four? Such as those in \cite[\S5]{CV}, double Ore extensions 
in~\cite{ZZ1,ZZ2}, generalized Laurent polynomial rings  
in~\cite{C.Goetz.Shelton}, etc.

10.
Does the line scheme of a generic quadratic regular algebra of global 
dimension four have a particular form? Perhaps a union of elliptic curves?
Or, perhaps it contains at least one elliptic curve?

11.
Suppose $A$ is as in Definition~\ref{ptmod} and $F$ is as in Theorem~\ref{ci}.
Let $I_k = \la f_1, \ldots , f_k\ra$ for all $k\leq n$, and let
$\widehat{\V(I_k)}$ denote the set of isomorphism classes of right base-point 
modules over $A/I_k$.
If $A$ is commutative, then, for each $k$, $\widehat{\V(I_k)}$ is a scheme, 
and so has a dimension.
In particular, if $A$ is the polynomial ring, then $F$ is regular if and
only if $\dim(\widehat{\V(I_k)}) = n- k - 1$, for all $k \leq n$. 
However, if $A$ is not commutative, is there an analogous statement and under 
what hypotheses on $A$ could it hold?

\medskip

\subsection{Related Topics}\hfill\\
\indent
Since the publication of \cite{ATV1}, the subject has branched out in many 
different directions, with the key topics being: classification of regular 
algebras; classification of projective surfaces; seeing which commutative 
techniques (e.g., blowing-up, blowing-down) carry over to the non-commutative
setting; and connections with differential geometry (e.g., via Poisson 
geometry).  Module categories and homological algebra provide a unifying 
umbrella over these topics. These different directions are highlighted 
in the references cited in the Introduction and throughout the text, and in 
the presentations from the 2013 MSRI program found in this journal issue.

New directions continue to emerge, with one of the most recent trends being 
the study of regular algebras and Hopf algebras together via the consideration
of Hopf actions on regular algebras, such as the work in \cite{CWWZ}. However, 
perhaps the most recent exciting triumph of the
subject is when the universal enveloping algebra of the Witt algebra was
viewed through the geometric lens of~\cite{ATV1} by Sierra and Walton, 
in~\cite{SW}, enabling them to solve the long-standing problem of whether
or not that algebra is noetherian.

In view of all these advances, it is now clear that the marriage of 
non-commutative algebra and algebraic geometry, \`a la~\cite{ATV1}, 
is a dynamic and evolving field of research.



\bigskip
\bigskip

\textheight8.85in

\end{document}